\documentclass{amsart}
\usepackage{amscd}

\theoremstyle{plain}
\newtheorem{theorem}{Theorem}[section]
\newtheorem{lemma}[theorem]{Lemma}
\newtheorem{proposition}[theorem]{Proposition}
\newtheorem{corollary}[theorem]{Corollary}
\theoremstyle{remark}
\newtheorem{remark}[theorem]{Remark}

\newcommand{\Z}{\mathbb{Z}}
\newcommand{\N}{\mathbb{N}}

\DeclareMathOperator{\GL}{GL}

\DeclareMathOperator{\tr}{tr}
\DeclareMathOperator{\Lie}{Lie}
\DeclareMathOperator{\Hom}{Hom}

\DeclareMathOperator{\Ker}{Ker}

\DeclareMathOperator{\Aut}{Aut}
\DeclareMathOperator{\Endom}{End}
\DeclareMathOperator{\Ext}{Ext}

\DeclareMathOperator{\Rep}{Rep}
\DeclareMathOperator{\Mat}{Mat}

\DeclareMathOperator{\depth}{depth}
\DeclareMathOperator{\codim}{codim}

\newcommand{\pr}{\Pi^\lambda}
\newcommand{\dbslash}{/\!\!/}
\DeclareMathOperator*{\rightoverleft}{\parbox{2em}{\centerline{$\longrightarrow$}\vskip -6pt\centerline{$\longleftarrow$}}}

\begin{document}
\title[Normality of Marsden-Weinstein reductions]{Normality of Marsden-Weinstein reductions for representations of quivers}

\author{William Crawley-Boevey}
\address{Department of Pure Mathematics, University of Leeds, Leeds LS2 9JT, UK}
\email{w.crawley-boevey@leeds.ac.uk}

\begin{abstract}
We prove that the Marsden-Weinstein reductions
for the moment map associated to representations
of a quiver are normal varieties. We give an application to
conjugacy classes of matrices.
\end{abstract}
\maketitle

\section{Introduction}
Throughout, we work over an algebraically closed field $K$ of characteristic zero.
Kraft and Procesi have proved in \cite{KP} that the closure of any conjugacy 
class of $n\times n$ matrices is a normal variety. They reduce to the case of
a nilpotent conjugacy class $C$, realize the closure $\overline{C}$ 
as an affine quotient variety $Z \dbslash H$, and then use Serre's criterion to prove that
$Z$ is normal. Here $Z$ is the variety of linear maps
\[
K^{n_1} 
\rightoverleft^{A_1}_{B_1}
K^{n_2} 
\rightoverleft^{A_2}_{B_2}
\cdots
\rightoverleft^{A_{t-2}}_{B_{t-2}} 
K^{n_{t-1}} 
\rightoverleft^{A_{t-1}}_{B_{t-1}} 
K^n 
\]
satisfying $B_1 A_1 = 0$ and $B_i A_i = A_{i-1} B_{i-1}$ for $2\le i \le t-1$, and
\[
H = \prod_{i=1}^{t-1} \GL(n_i,K)
\]
acting by conjugation (for suitable integers $t$ and $n_1 < n_2 < \dots < n_{t-1} < n$). 

In this paper we take these ideas very much further, and prove the normality 
of quotients similar to $Z \dbslash H$, but associated to any quiver. 

Let $Q$ be a quiver with vertex set $I$, and 
let $h(a)$ and $t(a)$ denote the head and tail 
vertices of an arrow $a\in Q$.
Let $\overline{Q}$ be the double of $Q$, 
obtained from $Q$ by adjoining a reverse arrow $a^*$ 
for each arrow $a$ in $Q$.
If $\alpha\in\N^I$, then the space of representations of $\overline{Q}$
of dimension vector $\alpha$ is
\[
\Rep(\overline{Q},\alpha) = 
\bigoplus_{a\in \overline{Q}}\Mat(\alpha_{h(a)}\times\alpha_{t(a)},K),
\]
the group
\[
\GL(\alpha) = \prod_{i\in I} \GL(\alpha_j,K)
\]
acts on it by conjugation,
its Lie algebra is
\[
\Endom(\alpha) = \prod_{i\in I} \Mat(\alpha_j,K),
\]
and there is a map
\[
\mu_\alpha:\Rep(\overline{Q},\alpha) \to \Endom(\alpha),
\quad
\mu_\alpha(x)_i = 
\sum_{\substack{a\in Q \\ h(a)=i}} x_a x_{a^*}
-
\sum_{\substack{a\in Q \\ t(a)=i}} x_{a^*} x_a,
\]
which is in fact a moment map for the action of 
$\GL(\alpha)$ on $\Rep(\overline{Q},\alpha)$ 
with respect to a natural symplectic form on $\Rep(\overline{Q},\alpha)$,
and after identifying the Lie algebra with its dual.

Identifying $\lambda\in K^I$ with the element of
$\Endom(\alpha)$ whose $i$th component is $\lambda_i$ times
the identity matrix of size $\alpha_i$, one can study the
fibre $\mu_\alpha^{-1}(\lambda)$ and the affine quotient
variety $N_Q(\lambda,\alpha) = \mu_\alpha^{-1}(\lambda) \dbslash \GL(\alpha)$.
Such a quotient is an appropriate analogue in the category of 
affine varieties of a `symplectic quotient', or 
`Marsden-Weinstein reduction' \cite{MW}. 
Note that this variety is independent of the orientation of $Q$.
(One might consider Marsden-Weinstein reductions associated to more 
general coadjoint orbits, but one obtains nothing new, see Remark \ref{r:notnew}.)

Constructions of this form were used by Kronheimer \cite{Kronheimer} 
to obtain the Kleinian singularities and their deformations and
resolutions (see also \cite{CS,CBH}), and they also play 
an important role in Nakajima's geometric construction 
of integrable representations of Kac-Moody algebras
\cite{NakajI}. For this reason we began a systematic study of
them in \cite{CBmm}. (See the introduction to that paper for a
detailed explanation of the connection with Nakajima's quiver varieties.)
Since then they have found applications back to Kleinian 
singularities \cite{CBex}, to the additive version of a problem of
Deligne and Simpson \cite{CBds}, and to Kac's conjectures on
absolutely indecomposable representations of quivers \cite{CBV}.
For some connections with noncommutative geometry, see \cite{BGK,BLB}.

Often the varieties are smooth, but this is not always the case, 
and it is an advantage of working in the 
algebraic geometric setting that one can deal with singular examples.
In general we expect to obtain an interesting class of singularities, often
having smooth deformations and resolutions 
with symplectic forms. (See the discussion in Section \ref{s:completion}.)
Our main result is as follows.

\begin{theorem}
\label{t:mainth}
$N_Q(\lambda,\alpha)$ is a normal variety.
\end{theorem}

This result does not appear to apply to the situation discussed at the start of 
the introduction, since the group $\GL(n,K)$ is missing from $H$, and
the relation $A_{t-1} B_{t-1}=0$ is also missing. However, these cancel out:
by separating the rows of $A_{t-1}$ into $n$ matrices of size $1\times n_{t-1}$, and the
columns of $B_{t-1}$ into $n$ matrices of size $n_{t-1}\times 1$, one can see that
$Z \cong \mu_\alpha^{-1}(0)$ and $Z\dbslash H \cong N_Q(0,\alpha)$
for a quiver $Q$ of shape
\[
\begin{picture}(170,20)
\put(0,10){\circle*{2.5}}
\put(30,10){\circle*{2.5}}
\put(63,10){\circle*{1}}
\put(70,10){\circle*{1}}
\put(77,10){\circle*{1}}
\put(110,10){\circle*{2.5}}
\put(140,10){\circle*{2.5}}
\put(170,10){\circle*{2.5}}
\put(5,10){\vector(1,0){20}}
\put(35,10){\vector(1,0){20}}
\put(85,10){\vector(1,0){20}}
\put(115,10){\vector(1,0){20}}
\put(145,16){\vector(1,0){20}}
\put(145,12){\vector(1,0){20}}
\put(148,8){\circle*{1}}
\put(153,8){\circle*{1}}
\put(158,8){\circle*{1}}
\put(145,4){\vector(1,0){20}}
\end{picture}
\]
with vertices $1,\dots,t$, with $n$ arrows between the last two, 
and $\alpha=(n_1,\dots,n_{t-1},1)$. (This method of `deframing' is 
discussed at the end of the introduction to \cite{CBmm}.)

We remark that unlike Kraft and Procesi's case, in the general case of 
Theorem~\ref{t:mainth} the fibre $\mu_\alpha^{-1}(\lambda)$ itself need not be 
normal; see the discussion at the end of \cite{CB} for an example.

Using Kraft and Procesi's construction of conjugacy class closures
several times, or Lemma \ref{l:kpgen} in the non-nilpotent case,
we have the following immediate consequence.

\begin{corollary}
If $C_1,\dots,C_k$ are conjugacy classes in $\Mat(n,K)$, then the
quotient variety
\[
\left\{ (M_i) \in \prod_{i=1}^k \overline{C_i} \: \Bigg\vert \: \sum_{i=1}^k M_i=0 \right\} \dbslash \GL(n,K)
\]
is normal.
\end{corollary}

Indeed the quotient is isomorphic to a variety of the form $N_Q(\lambda,\alpha)$ where $Q$ is a star with 
$k$ arms and $\alpha$ is a dimension vector which is $n$ at the central 
vertex. The vanishing of the component of the moment map for the central vertex 
corresponds to the relation $\sum_{i=1}^k M_i=0$.

\medskip
In order to describe our earlier work on $N_Q(\lambda,\alpha)$, and for later use, 
we introduce some notation.
Associated to a quiver $Q$ with vertex set $I$ there is a quadratic form
\[
q(\alpha) = \sum_{i\in I} \alpha_i^2 - \sum_{a\in Q} \alpha_{h(a)}\alpha_{t(a)}
\]
on $\Z^I$. 
We denote by $(\alpha,\beta)$ the symmetric bilinear form with $(\alpha,\alpha)=2q(\alpha)$,
and we define $p(\alpha)=1-q(\alpha)$.
(If we need to specify the quiver, we decorate $q$, $(-,-)$ and $p$ with
subscripts.)
Associated to $Q$ there is a subset of $\Z^I$ of `roots' for $Q$ (see for example \cite{Kac}), 
and if $\lambda\in K^I$
we denote by $R_\lambda^+$ the set of positive roots $\alpha$ with 
$\lambda\cdot\alpha := \sum_{i\in I}\lambda_i \alpha_i = 0$. 
We define $\Sigma_\lambda$ to be the set of $\alpha\in R_\lambda^+$
with the property that 
\[
p(\alpha) > p(\beta^{(1)}) + p(\beta^{(2)}) + \dots
\]
whenever $\alpha$ can be written as a sum $\alpha=\beta^{(1)} + \beta^{(2)} + \dots$
of two or more elements of~$R_\lambda^+$.

The elements of $\mu_\alpha^{-1}(\lambda)$ correspond to representations 
of dimension vector $\alpha$ of a certain algebra, the 
\emph{deformed preprojective algebra} $\pr$ of \cite{CBH}, and the quotient 
$N_Q(\lambda,\alpha)$ classifies the isomorphism classes of semisimple 
representations. We showed in \cite{CBmm} that there is a simple representation 
of $\pr$ of dimension vector $\alpha$ if and only if $\alpha\in\Sigma_\lambda$.
Moreover, in this case $\mu_\alpha^{-1}(\lambda)$ and $N_Q(\lambda,\alpha)$
are irreducible varieties. 
In \cite{CBde} we studied $N_Q(\lambda,\alpha)$ for arbitrary $\alpha$,
and showed that $N_Q(\lambda,\alpha)$ is either an irreducible variety
or empty, according to whether or not $\alpha$ can be written as a 
sum of elements of~$R_\lambda^+$.

Note that, as in \cite{CBde}, we consider $N_Q(\lambda,\alpha)$ as a variety,
that is, we endow it with the reduced scheme structure. In \cite{CBmm} we worked 
with schemes, and showed that the natural scheme structures on 
$\mu_\alpha^{-1}(\lambda)$ and $N_Q(\lambda,\alpha)$
are reduced when $\alpha\in\Sigma_\lambda$. In general, however, this is
not clear.

\medskip
The contents of this paper are as follows. In Section \ref{s:sympbimod} we study 
symplectic forms on bimodules over a semisimple algebra, and in 
Theorem \ref{t:maxisotrop} we prove an analogue of the (easy, linear) Darboux Theorem.
In Section \ref{s:coordfree} we use this structure theory to give a 
coordinate-free construction of the varieties $N_Q(\lambda,\alpha)$.
This is of interest in its own right, since this is how the varieties often arise,
but our main application of it comes in 
Section \ref{s:etalelocal}, when we study the etale local structure of 
$N_Q(\lambda,\alpha)$. In Corollary \ref{c:etlocal}
we show that any point in one of the varieties $N_Q(\lambda,\alpha)$ 
has a neighbourhood which is isomorphic in the etale topology to a neighbourhood 
of the origin in $N_{Q'}(0,\kappa)$ for some $Q'$ and~$\kappa$.

There is a stratification of $N_Q(\lambda,\alpha)$ according to `representation type' 
\cite{CBmm}
with the stratum of type $(k_1,\beta^{(1)}; \dots; k_r, \beta^{(r)})$ 
containing those semisimple representations of $\pr$ which involve $r$ 
non-isomorphic simples, with dimension vectors $\beta^{(t)}$, and multiplicities $k_t$. 
An important problem is to compute the dimension
of the inverse image of this stratum in $\mu_\alpha^{-1}(\lambda)$.
In Section \ref{s:modexts} we prove some results 
which work for an arbitrary finitely generated algebra, and then in 
Section \ref{s:dimstrata} we apply this to the preprojective algebra $\Pi^0$
to get a bound on the dimensions of inverse images of strata in case $\lambda=0$.
This partially generalizes Lusztig's result \cite{Lusztig} that the nilpotent cone
in $\mu_\alpha^{-1}(0)$ is Lagrangian in $\Rep(\overline{Q},\alpha)$
when $Q$ has no oriented cycles.

In Section \ref{s:normquot} we prove a normality condition for quotient
schemes $X\dbslash G$ with $G$ a reductive group acting on a scheme $X$.
If $X$ is normal, then it is well-known that $X\dbslash G$ is normal.
Recall that Serre's criterion says that a scheme is normal if and only if
it has the regularity condition $(R_1)$ and the depth condition $(S_2)$.
In our situation, $X$ has $(S_2)$, but may fail $(R_1)$. However $X\dbslash G$ has
$(R_1)$. Can one deduce that $X\dbslash G$ is normal? 
In Corollary \ref{c:quotnormal} we show that the answer is `yes' if
$X \dbslash G$ has an open subscheme $U$ which is normal,
whose complement $Z$ has codimension at least two, 
and such that the inverse image of $Z$ in $X$ has
codimension at least two.

In Section \ref{s:completion} we complete the proof of
our main theorem. It is amusing that the cases that cannot be dealt with
by the etale local structure result and an induction, are nearly Kleinian,
and hence trivially normal. Finally we show that in many
cases the varieties $N_Q(0,\alpha)$ are genuinely singular.

For further developments based on the results in this paper we refer to \cite{LB1,LB2}.

I would like to thank M. P. Holland for his assistance and encouragement during 
this work. Much of the research was done in Spring 2000 while visiting first 
the program on `Noncommutative Algebra' at MSRI, and then the Sonderforschungsbereich 
on `Discrete Structures in Mathematics' at Bielefeld University. I would like to thank 
my hosts at both institutions for their hospitality.

\section{Symplectic forms on bimodules}
\label{s:sympbimod}

Let $A$ be a finite-dimensional semisimple associative $K$-algebra. Since
$K$ is algebraically closed, $A$ is isomorphic
to a direct sum of full matrix algebras. It follows that any $A$-$A$-bimodule
is semisimple, so any sub-bimodule of an $A$-$A$-bimodule has
a complement which is a sub-bimodule.

If $M$ is a finite-dimensional $A$-$A$-bimodule, then a 
symplectic bilinear form $\omega$ on $M$ is said to be
\emph{balanced} if $\omega(x,ay) = \omega(xa,y)$
for all $a\in A$ and $x,y\in M$. Using the skew-symmetry
of $\omega$ one has $\omega(ax,y) = -\omega(y,ax) = -\omega(ya,x) = \omega(x,ya)$.

If $U$ is a subspace of $M$, one defines 
$U^\perp=\{ m\in M \mid \omega(m,u)=0 \text{ for all $u\in U$}\}$. 
The non-degeneracy of $\omega$ implies 
that $\dim U^\perp = \dim M - \dim U$ and $U^{\perp\perp} = U$.
If $U$ is a sub-bimodule then so is $U^\perp$.
Recall that $U$ is \emph{isotropic} (respectively \emph{coisotropic}, 
respectively \emph{Lagrangian}) if $U \subseteq U^\perp$
(respectively $U \supseteq U^\perp$, respectively $U = U^\perp$).

\begin{lemma}\label{l:simpisotr}
In an $A$-$A$-bimodule with a balanced symplectic form, any
simple sub-bimodule must be isotropic.
\end{lemma}

\begin{proof}
Suppose that $T$ is a simple sub-bimodule of $M$, and that $M$ has
a balanced symplectic form $\omega$. Then either $T\cap T^\perp = T$,
or $T\cap T^\perp = 0$. In the former case $T$ is isotropic, so
suppose the latter holds, that is, that the restriction of $\omega$
to $T$ is non-degenerate, hence a symplectic form on $T$.

Since $A$ is a product of full matrix algebras, $T$ must be isomorphic
to a bimodule of the form $X\otimes_K Y^*$
where $X$ and $Y$ are simple left $A$-modules.
The symplectic form gives a bimodule 
isomorphism $T\to T^*$, so $X\cong Y$.
It follows that $T \cong \Endom_K(X)$.
Since $X$ is simple, the homothety $A\to \Endom_K(X)$
is surjective, and since this is a bimodule map, $T$ is a quotient of $A$ as a bimodule.

Let $t$ be the image in $T$ of the identity element of $A$. Then
every element of $T$ is of the form $at=ta$ for some $a\in A$.
Now $\omega(at,t) = \omega(ta,t) = \omega(t,at) = -\omega(at,t)$.
Thus $\omega(at,t)=0$, and since this is true for all $a$,
non-degeneracy implies that $t=0$, a contradiction.
\end{proof}

We now have our analogue of the Darboux Theorem.

\begin{theorem}\label{t:maxisotrop}
If $M$ is an $A$-$A$-bimodule with a balanced 
symplectic form $\omega$, then any maximal isotropic sub-bimodule $S$ 
of $M$ has an isotropic bimodule complement. 
Both $S$ and the complement are Lagrangian.
Moreover there is an isomorphism of bimodules $M\cong S\oplus S^*$, 
under which the symplectic form is given by
$\omega((s,f),(s',f')) = f(s')-f'(s)$.
\end{theorem}

\begin{proof}
Since bimodules over $A$ are semisimple, $S$ has a bimodule complement $C$. 
The fact that $\omega$ is balanced implies that the map $C\to S^*$ induced by $\omega$
is a bimodule map. It is injective, for if $T$ were a simple sub-bimodule 
of the kernel, then $S\oplus T$ would be isotropic by Lemma \ref{l:simpisotr}, 
contradicting maximality.

Dualizing, the map $f:S\to C^*$ induced by $\omega$ is surjective. 
Semisimplicity for the category of bimodules then implies that
the map $C\to C^*$ induced by $\omega$ factors through $f$
by a bimodule map, say $\theta:C\to S$. Thus $\omega(\theta(c),c') = \omega(c,c')$
for all $c,c'\in C$. By skew-symmetry, also $\omega(c',\theta(c)) = \omega(c',c)$.

Let $D = \{ c-\frac12 \theta(c) \mid c\in C\}$. Clearly this is a bimodule
complement to $S$. Moreover
\[
\omega(c - \frac12 \theta(c), c' - \frac12 \theta(c'))
= \omega(c, c')
- \omega(\frac12 \theta(c), c')
- \omega(c, \frac12 \theta(c'))
+ \omega(\frac12 \theta(c),\frac12 \theta(c')).
\]
Now the last term on the right hand side is zero since $S$ is isotropic, and
the second and third terms each give $-\frac12 \omega(c,c')$.
Thus the right hand side vanishes.

Now $S$ and $D$ are isotropic, so both have dimension bounded by $\frac12 \dim M$.
But this implies that they have dimension equal to this, and hence both are
Lagrangian. Now the map $D\to S^*$ induced by $\omega$ is an isomorphism, and
this gives the final assertion.
\end{proof}

\begin{corollary}\label{c:isotcompl}
If $M$ is an $A$-$A$-bimodule with a balanced 
symplectic form $\omega$, then any isotropic sub-bimodule $S$ of $M$ has a 
coisotropic bimodule complement.
\end{corollary}

\begin{proof}
Extend $S$ to a maximal isotropic sub-bimodule, which by semisimplicity we may write
as $S\oplus T$, and choose a Lagrangian complement $C$ to this. Then $T\oplus C$
is a complement to $S$, and it is coisotropic since
$(T\oplus C)^\perp \subseteq C^\perp = C \subseteq T\oplus C$.
\end{proof}

\section{A coordinate-free description of $N_Q(\lambda,\alpha)$}
\label{s:coordfree}

We consider quadruples $(A,M,\omega,\zeta)$ where 
$A$ is a finite-dimensional semisimple associative $K$-algebra,
$M$ is a finite-dimensional $A$-$A$-bimodule,
$\omega$ is a balanced symplectic bilinear form on $M$, and $\zeta:A\to K$ is a linear map
with $\zeta(ab) = \zeta(ba)$ for all $a,b\in A$ (a `trace function').

Let $G$ be the group of units of $A$. By Wedderburn's Theorem 
it is a product of copies of the general linear group. 
It acts by conjugation on any $A$-$A$-bimodule $N$ via $g\cdot n = gng^{-1}$ 
(where the multiplication on the right hand side of the equation 
is the bimodule action of $A$ on $N$).
Recall that the \emph{centre} of $N$
is the set 
\[
\{ n\in N \mid \text{$an = na$ for all $a\in A$}\}.
\]
The fact that $G$ is Zariski-dense in $A$ implies that the centre of $N$ 
coincides with the set of $G$-invariant elements of $N$.
In particular, the condition imposed on $\zeta$ is equivalent to insisting
that it be a $G$-invariant element of $A^*$.

Let $\mu:M\to A^*$ be the map defined by $\mu(m)(a) = \omega(m,am)$.
Clearly we have $\omega(m,am) = \omega(ma,m) = -\omega(m,ma)$,
so that $\mu(m)(a) = \frac{1}{2} \omega(m,[a,m])$ where $[a,m] =am-ma$.
Associated to the quadruple $(A,M,\omega,\zeta)$ there
is an affine variety
\[
N(A,M,\omega,\zeta) = \mu^{-1}(\zeta) \dbslash G.
\]
It is a Marsden-Weinstein reduction because of the following 
observation.

\begin{lemma}
$\mu$ is a moment map for the action of $G$ on $M$.
\end{lemma}

\begin{proof}
This means that $d\mu_m(v)(a) = \omega_m (v,a_m)$
where $m\in M$, $v\in T_m M$ and $a\in\Lie(G)$.
Here $\omega_m$ is the symplectic form induced by
$\omega$ on $T_m M$ and $m\mapsto a_m$ is the vector field
on $M$ induced by $a$.

Each tangent space of $M$ can be identified with
$M$ itself, and $\Lie(G)$ can be identified
with $A$. The vector field induced by $a$ is then the map
$M\to M$ sending $m$ to $[a,m]=am-ma$. The equation
thus becomes $d\mu_m (v)(a) = \omega(v,[a,m])$ for 
$m,v\in M$ and $a\in A$.

Now the formula for $\mu$ immediately gives
$d\mu_m (v)(a) = \omega(v,am)+\omega(m,av)$,
and using the properties of $\omega$ this becomes
$\omega(v,[a,m])$, as required.
\end{proof}

A quiver $Q$ with vertex set $I$, and elements $\alpha\in\N^I$ and
$\lambda\in K^I$ give rise to a quadruple
$(\Endom(\alpha),\Rep(Q,\alpha),\omega_\alpha,\zeta_\lambda)$
where $\omega_\alpha$ is the canonical symplectic
form on $\Rep(Q,\alpha)$,
\[
\omega_\alpha(x,y) = \sum_{a\in Q} \bigl(\tr(x_{a^*} y_a) - \tr(x_a y_{a^*}) \bigr),
\]
and $\zeta_\lambda$ is defined by 
$\zeta_\lambda(\theta) = \sum_{i\in I} \lambda_i \tr (\theta_i)$.

\begin{lemma}
$N(\Endom(\alpha),\Rep(Q,\alpha),\omega_\alpha,\zeta_\lambda) \cong N_Q(\lambda,\alpha)$.
\end{lemma}

\begin{proof}
Using the trace pairing 
$(\theta,\phi) \mapsto \sum_{i\in I} \tr(\theta_i \phi_i)$
to identify $\Endom(\alpha)$ with its dual, the map 
$\mu$ coincides with the map $\mu_\alpha$ and $\zeta_\lambda$ 
coincides with $\lambda$ as an element of $\Endom(\alpha)$.
\end{proof}

We say that two quadruples $(A,M,\omega,\zeta)$
and $(A',M',\omega',\zeta')$ are \emph{equivalent} if
there is an algebra isomorphism $f:A\to A'$ 
with $\zeta'(f(a)) = \zeta(a)$ for all $a\in A$,
and an $A$-$A$-bimodule isomorphism $g:M\to {}_f {M'}_{f}$
with $\omega'(g(m_1),g(m_2)) = \omega(m_1,m_2)$
for all $m_1,m_2\in M$.
Clearly an equivalence of quadruples gives an isomorphism 
$N(A,M,\omega,\zeta)\cong N(A',M',\omega',\zeta')$.
Our coordinate-free description of the spaces $N_Q(\lambda,\alpha)$ 
is completed with the following result.

\begin{lemma}
\label{l:quadquiv}
Every quadruple is equivalent to 
$(\Endom(\alpha),\Rep(Q,\alpha),\omega_\alpha,\zeta_\lambda)$
for some $Q$, $\alpha$ and $\lambda$.
\end{lemma}

\begin{proof}
Let $(A,M,\omega,\zeta)$ be a quadruple.
Choose a maximal isotropic sub-bimodule $S$ of $M$
and identify $M$ with $S\oplus S^*$ as
in Theorem \ref{t:maxisotrop}.

The vertex set $I$ of $Q$ is chosen to 
index the simple left $A$-modules
up to isomorphism. The dimension vector
$\alpha$ is defined by letting $\alpha_i$ be the dimension
of the $i$th simple module. Clearly $A\cong \Endom(\alpha)$.

If $X$ and $Y$ are simple $A$-modules, then $X\otimes_K Y^*$
is a simple $A$-$A$-bimodule, and all simple modules are obtained in
this way. We define the quiver $Q$ by letting the number of arrows 
from $X$ to $Y$ be the multiplicity of this simple module in $S$.

Finally, $\zeta$ is a $G$-invariant element of $A^* \cong \Endom(\alpha)^*$. Using
the trace pairing it corresponds to a $G$-invariant element
of $\Endom(\alpha)$. The $G$-invariance means that each of the matrices
involved is a scalar multiple of the identity matrix, and combining
the scalars gives the element $\lambda\in K^I$with $\zeta_\lambda=\zeta$.
\end{proof}

\section{Etale local structure}
\label{s:etalelocal}

In this section we study the local structure of our 
Marsden-Weinstein reductions. In a differential geometric setting
this has already been done in a special case by Kronheimer \cite[Lemma 3.3]{Kronheimer}.
In our situation we have the added complication of etale morphisms,
so need Luna's Slice Theorem \cite{Luna}.
(See also \cite[\S 6]{NakajI} and \cite{LBP}.)
If $G$ is a reductive group acting on an affine variety $X$,
we use the notation $X\dbslash G$ where Luna writes $X/G$, but otherwise
we use Luna's notation. We denote by $\pi_X$ the natural map $X\to X\dbslash G$.

Let $(A,M,\omega,\zeta)$ be a quadruple as in the previous section.
Let $G$ be the group of units of $A$.
Let $\mu:M\to A^*$ be the moment map $\mu(m)(a) = \omega(m,am)$.
Fix $x\in\mu^{-1}(\zeta)$ with the property that its $G$-orbit is closed.

\begin{lemma}
$[A,x]$ is an isotropic subspace of $M$. That is, $[A,x] \subseteq [A,x]^\perp$.
\end{lemma}

\begin{proof}
Since $\omega$ is balanced and skew-symmetric,
for $a,b\in A$ we have
\begin{align*}
\omega([a,x],[b,x]) 
&= \omega(ax,bx) - \omega(ax,xb) - \omega(xa,bx) + \omega(xa,xb) \\
&= \omega(axb,x) + \omega(xb,ax) - \omega(x,abx) + \omega(x,axb) \\
&= \omega(x,bax) - \omega(x,abx).
\end{align*}
But this can be rewritten as $\mu(x)(ba) - \mu(x)(ab) = \zeta(ba) - \zeta(ab) = 0$.
\end{proof}

\begin{lemma}
$A_x = \{a\in A \mid ax = xa \}$ is a semisimple subalgebra of $A$.
\end{lemma}

\begin{proof}
The stabilizer $G_x$ of $x$
is reductive by Matsushima's Theorem, see \cite{Luna}.
Now $G_x$ is the group of units of $A_x$, and
the elements $1+a$ with $a$ in the radical 
of $A_x$ form a closed connected unipotent 
normal subgroup. Thus the radical must be zero.
\end{proof}

Choose an $A_x$-$A_x$-bimodule complement
$L$ to $A_x$ in $A$, so that $A = A_x\oplus L$.
Observe that $[A,x]$ is an $A_x$-$A_x$-sub-bimodule of $M$.
By Corollary \ref{c:isotcompl} we can 
choose a coisotropic $A_x$-$A_x$-bimodule complement 
$C$ to $[A,x]$ in $M$, so $M = [A,x]\oplus C$.
Let $W = C\cap [A,x]^\perp$.
Let $\mu_x:M\to A_x^*$ be the map obtained by composing
$\mu$ with the restriction map $A^*\to A_x^*$.
Let $\hat{\mu}$ be the restriction of $\mu_x$ to $W$.

\begin{lemma}
The restriction of the symplectic form $\omega$ to $W$ is non-degenerate, 
and hence a symplectic form. It is balanced for the 
$A_x$-$A_x$-bimodule structure, and
the corresponding moment map is $\hat{\mu}$.
\end{lemma}

\begin{proof}
Let $0\neq h\in W$. We show that
there is an element $h'\in W$ with $\omega(h,h')\neq 0$.
Since $\omega$ is non-degenerate, $C^\perp \cap [A,x]^\perp = 0$,
so by dimensions $M = C^\perp \oplus [A,x]^\perp$.
Non-degeneracy now implies that there 
are elements $p\in C^\perp$ and $q\in [A,x]^\perp$
with $\omega(h,p+q)\neq 0$. Thus $\omega(h,q)\neq 0$.
Now since $M = [A,x]\oplus C$ and $[A,x]\subseteq [A,x]^\perp$ 
we have $[A,x]^\perp = [A,x] \oplus W$.
The claim follows on writing $q = q'+h'$ corresponding to 
this decomposition, and observing that $\omega(h,q')=0$. 
The rest is clear.
\end{proof}

Define $\nu:C\to L^*$ by 
$\nu(c)(\ell) = \omega(c,\ell x) + \omega(c,\ell c) + \omega(x,\ell c)$. 
Thus
\[
\mu(x+c)(a+\ell)=\zeta(a+\ell)+\mu_x(c)(a)+\nu(c)(\ell)
\]
for $c\in C$, $a\in A_x$ and $\ell\in L$.

\begin{lemma}\label{l:et1}
The assignment $c\mapsto x+c$ induces a $G_x$-equivariant map
$\mu_x^{-1}(0)\cap \nu^{-1}(0) \to \mu^{-1}(\zeta)$,
and the induced map
$(\mu_x^{-1}(0)\cap \nu^{-1}(0)) \dbslash G_x \to \mu^{-1}(\zeta) \dbslash G$
is etale at $0$.
\end{lemma}

\begin{proof}
We consider the map $C\to M$ sending $c$ to $x+c$. 
The formula for $\mu(x+c)$ shows that 
$x+c\in\mu^{-1}(\zeta)$ if and only if $c\in\mu_x^{-1}(0)\cap \nu^{-1}(0)$,
from which the first assertion follows.
The conjugation action of $G$ on $M$ gives a map $G\times C\to M$,
and the fact that $M = [A,x]\oplus C$ implies that the differential 
of this map at $(1,0)$ is surjective. Let $X = G\times_{G_x} C$.
(See \cite[p86]{Luna}.) By dimensions it follows that the induced map
$\phi:X\to M$ is etale at $\overline{(1,0)}$.

Thus Luna's Fundamental Lemma \cite[p94]{Luna} applies to $\phi$.
There is an affine open subset $U$ of $X$ which 
contains $\overline{(1,0)}$, is saturated for $\pi_X$, 
the restriction of $\phi$ to $U$ is etale,
the image $V = \phi(U)$ is an affine open subset of $M$,
saturated for $\pi_M$, 
the morphism $\phi / G : U \dbslash G\to V \dbslash G$ is etale,
and the morphisms $\phi:U\to V$ and $\pi_U$
induce a $G$-isomorphism $U\to V\times_{V\dbslash G} U\dbslash G$.

Let $V' = V\cap\mu^{-1}(\zeta)$. Since $V'$ is a closed subvariety 
of $V$, one can identify $U' = V'\times_V U$ with a closed subvariety 
of $U$, and clearly $U' = U \cap X'$
where $X' = G \times_{G_x} C'$ and $C' = \mu_x^{-1}(0)\cap \nu^{-1}(0)$.

Now \cite[Lemma, p95]{Luna}, applied to the map $\phi:U\to V$
and the subvariety $V'$, implies that 
the morphism 
\[
\phi'/G : (U \cap X')\dbslash G \to (V\cap\mu^{-1}(\zeta))\dbslash G
\]
is etale. 

The fact that $U$ is saturated for $\pi_X$ implies that $U \cap X'$
is saturated for $\pi_{X'}$, so the domain of $\phi'/G$ is an open
subset of $X'\dbslash G \cong C' \dbslash G_x$. Similarly,
the fact that $V$ is saturated for $\pi_M$ implies that $V\cap\mu^{-1}(\zeta)$
is saturated for $\pi_{\mu^{-1}(\zeta)}$, so the codomain of $\phi'/G$ is an open
subset of $\mu^{-1}(\zeta)\dbslash G$.
The claim follows.
\end{proof}

\begin{lemma}
$\dim W = \dim C - \dim L$.
\end{lemma}

\begin{proof}
Since $[A,x] \subseteq [A,x]^\perp$ we have $C + [A,x]^\perp = M$. Thus
\[
\dim W = \dim C + \dim [A,x]^\perp - \dim M = \dim C - \dim L,
\]
where the last equality
holds because $\dim [A,x]^\perp = \dim M - \dim [A,x] = \dim M - \dim A + \dim A_x = \dim M - \dim L$.
\end{proof}

\begin{lemma}\label{l:nusmooth}
The map $\nu:C\to L^*$ is smooth at 0, so $\nu^{-1}(0)$ is smooth at 0. Its tangent space at
0 is identified with $W$.
\end{lemma}

\begin{proof}
Identifying $C$ and $L^*$ with their tangent spaces 
at the origin, the differential of $\nu$ at the origin
is given by
\[
d\nu_0(c)(\ell) = \omega(c,\ell x)+\omega(x,\ell c) = \omega(c,[\ell,x])
\]
so $\Ker d\nu_0 = C\cap [L,x]^\perp = W$.
By the dimension formula $d\nu_0$ must be surjective, 
so $\nu$ is smooth at 0. 
\end{proof}

\begin{lemma}\label{l:cdecomp} 
$C = C^\perp \oplus W$
and 
$C^\perp = C \cap W^\perp$.
\end{lemma}

\begin{proof}
$W = C\cap [A,x]^\perp = C^{\perp\perp} \cap [A,x]^\perp = (C^\perp + [A,x])^\perp$,
so $W^\perp = C^\perp + [A,x]$.
Since $C$ is coisotropic, this implies that
$C\cap W^\perp = C^\perp + (C\cap [A,x]) = C^\perp$.
The left hand side can be rewritten as
$C^{\perp\perp}\cap W^\perp = (C^\perp+W)^\perp$,
and hence $C^\perp + W = C$.
Finally, $C^\perp\cap W  = C^\perp \cap [A,x]^\perp  = (C+[A,x])^\perp = M^\perp = 0$.
\end{proof}

\begin{lemma}\label{l:et2} 
There is a morphism 
$(\mu_x^{-1}(0)\cap \nu^{-1}(0)) \dbslash G_x\to \hat{\mu}^{-1}(0) \dbslash G_x$, 
sending 0 to 0, which is etale at 0.
\end{lemma}

\begin{proof}
Let $p:C\to W$ be the projection coming from
the decomposition of Lemma \ref{l:cdecomp}. 
Since $C^\perp$ is isotropic, it follows that $\hat{\mu} p = \mu_x$.
Let $\phi:\nu^{-1}(0)\to W$ be the restriction of $p$. 
By Lemma \ref{l:nusmooth} it is etale at 0.

We apply \cite[Lemma 3, p93]{Luna} to the $G_x$-morphism $\phi$.
Since the codomain is a vector space it is smooth,
so certainly normal. Thus there is an affine open subvariety $U'$ of
$\nu^{-1}(0)\dbslash G_x$ containing 0 such that the restriction
of 
\[
\phi/G_x : \nu^{-1}(0)\dbslash G_x \to W\dbslash G_x
\]
to $U'$ is etale. In particular $\phi/G_x$ is etale at 0.

Now since $\hat{\mu} p = \mu_x$ we have 
$\phi^{-1}(\hat{\mu}^{-1}(0)) = \mu_x^{-1}(0)\cap \nu^{-1}(0)$. 
Thus the pullback of $\phi / G_x$ along
the closed immersion $\hat{\mu}^{-1}(0) \dbslash G_x\to W\dbslash G_x$
is a morphism as in the statement of the lemma.
The assertion follows.
\end{proof}

Given points $x_i$ in varieties $X_i$ ($i=1,2$),
we say that the $x_i$ have neighbourhoods which are 
\emph{isomorphic in the etale topology} 
if there is a variety $Z$, a point $z\in Z$, and morphisms
from $Z$ to the $X_i$, etale at $z$ and sending $z$ to $x_i$.
For later use, note that $X_1$ is normal at $x_1$ if and only
if $X_2$ is normal at $x_2$ by \cite[I, Theorem 9.5]{SGA1}.

\begin{theorem}\label{th:ettopcoordfree}
The points $x\in \mu^{-1}(\zeta) \dbslash G$
and $0 \in \hat{\mu}^{-1}(0) \dbslash G_x$
have neighbourhoods which are isomorphic in the etale topology.
\end{theorem}

\begin{proof}
Combine Lemmas \ref{l:et1} and \ref{l:et2}.
\end{proof}

Combining this result with the coordinate-free description
of $N_Q(\lambda,\alpha)$ we obtain the following.

\begin{corollary}
\label{c:etlocal}
Suppose $Q$ is a quiver with vertex set $I$, $\lambda\in K^I$ and $\alpha\in \N^I$.
If $x\in N_Q(\lambda,\alpha)$, then there is a quiver $Q'$ and dimension vector $\kappa$ such
that the points $x\in N_Q(\lambda,\alpha)$ and $0\in N_{Q'}(0,\kappa)$
have neighbourhoods which are isomorphic in the etale topology.
\end{corollary}

Here $A = \bigoplus_{i\in I} \Mat(\alpha_i,K)$.
If $x$ has representation type $(k_1,\beta^{(1)};\dots ;k_r,\beta^{(r)})$, 
then 
\[
A_x \cong \Mat(k_1,K)\oplus\dots\oplus \Mat(k_r,K).
\]
Thus $Q'$ has vertex set $\{1,\dots,r\}$ and 
the dimension vector is $\kappa = (k_1,\dots,k_r)$.
To determine the arrows in $Q'$, or rather $\overline{Q'}$, one
needs to determine the multiplicities of the simple
$A_x$-$A_x$-bimodules $K^{k_s}\otimes_K (K^{k_t})^*$ in $W$. 
Now the natural map $W\to [A,x]^\perp / [A,x]$ 
is an isomorphism of $A_x$-$A_x$-bimodules,
and we have $[A,x] \cong A / A_x$ and $M/ [A,x]^\perp \cong (A/A_x)^*$. 
It thus suffices to determine how the simple $A$-$A$-bimodules
$K^{\alpha_i}\otimes_K (K^{\alpha_j})^*$ decompose into 
simple $A_x$-$A_x$-bimodules. This follows from the decomposition
\begin{align*}
K^{\alpha_i} 
&\cong (K^{\beta^{(1)}_i})^{k_1} \oplus \dots \oplus (K^{\beta^{(r)}_i})^{k_r} \\
&\cong (K^{k_1})^{\beta^{(1)}_i} \oplus \dots \oplus (K^{k_r})^{\beta^{(r)}_i}
\end{align*}
It follows that $\overline{Q'}$ has $2p_Q(\beta^{(i)})$ loops at vertex $i$, 
and $-(\beta^{(i)},\beta^{(j)})_Q$ arrows from $i$ to $j$, if $i\neq j$.

\section{Varieties of module extensions}
\label{s:modexts}

We prove a result about varieties of modules which is used in the next section.

Let $A$ be a finitely generated associative $K$-algebra,
let $I$ be a finite set, and let $e_i$ ($i\in I$) be a complete
set of orthogonal idempotents in $A$, meaning that
$e_i e_j = \delta_{ij} e_i$ and $\sum_{i\in I} e_i = 1$.
If $M$ is a finite dimensional $A$-module,
its \emph{dimension vector} is the
element of $\N^I$ whose $i$th component
is $\dim e_i M$.
For $\alpha\in\N^I$ we define $\Rep(A,\alpha)$
to be the variety of $A$-module structures on
$K^\alpha = \bigoplus_{i\in I} K^{\alpha_i}$ for which
$e_i$ acts as the projection onto the $i$th summand.
(Any such $A$ can be realized as a quotient of a path 
algebra $KQ$, with the $e_i$ corresponding to the trivial paths, 
and then $\Rep(A,\alpha)$ is a closed subvariety of $\Rep(Q,\alpha)$.)
If $x\in\Rep(A,\alpha)$ we denote the corresponding $A$-module by $K_x$.
The group
\[
\GL(\alpha) = \prod_{i\in I}\GL(\alpha_i,K)
\]
acts by conjugation on $\Rep(A,\alpha)$ and its orbits correspond to
isomorphism classes of $A$-modules $M$ of dimension vector $\alpha$.

\begin{lemma}\label{l:extE}
Suppose that $M$ and $L$ are finite-dimensional
$A$-modules and that $M\oplus L$ has dimension vector
$\alpha$. The set $\mathcal{E}$ of triples
\[
(x,\theta,\phi)\in \Rep(A,\alpha)\times \Hom_K(M,K^\alpha)\times \Hom_K(K^\alpha,L)
\]
such that $\theta$ and $\phi$ are $A$-module maps
for the module structure on $K^\alpha$ given by $x$, and 
the sequence
\[
0\to M \xrightarrow{\theta} K_x \xrightarrow{\phi} L\to 0
\]
is exact, is locally closed, and
$\dim\mathcal{E}
= \alpha\cdot\alpha + \dim\Ext^1(L,M) - \dim\Hom(L,M)$.
\end{lemma}

\begin{proof}
$\mathcal{E}$ is locally closed since the conditions that $\theta$ and $\phi$ 
be $A$-module maps and that $\phi\theta=0$ are closed, while the conditions 
that $\theta$ be injective and $\phi$ surjective are open.

Since $A$ is finitely generated as a $K$-algebra, it
follows that $L$ is a finitely-presented $A$-module,
so that $\Ext^1 (L,M)$ is finite dimensional.
Choose an exact sequence
\[
\xi:0\to U\xrightarrow{p} V\xrightarrow{q} L\to 0
\]
with the property that the connecting map 
$\Hom(U,M)\to \Ext^1(L,M)$ is surjective.
For example one can take $\xi$ to be the universal
exact sequence, with $U=M^k$ where $k=\dim\Ext^1(L,M)$.

Let $\mathcal{D}$ be the set of tuples 
\[
((x,\theta,\phi),f,g)
\in
\mathcal{E}\times\Hom(U,M)\times\Hom_K(V,K^\alpha)
\]
where $g$ is an $A$-module map for the module structure given by $x$, 
and the diagram
\[
\begin{CD}
0 @>>> U @>p>> V @>q>> L @>>> 0 \\
@. @VfVV @VgVV @\vert @. \\
0 @>>> M @>\theta>> K_x @>\phi>> L @>>> 0.
\end{CD}
\]
commutes. Clearly $\mathcal{D}$ is a closed subset.

We relate the dimensions of $\mathcal{E}$ and 
$\mathcal{D}$ by computing the fibres of the 
projection $\pi:\mathcal{D}\to\mathcal{E}$.
Suppose $(x,\theta,\phi)\in\mathcal{E}$.
By the surjectivity of the connecting map,
there is some $((x,\theta,\phi),f,g)\in\mathcal{D}$.
Moreover, given homomorphisms $f'\in\Hom(U,M)$ and 
$g'\in\Hom(V,K_x)$, we have
$((x,\theta,\phi),f+f',g+g')\in\mathcal{D}$
if and only if 
$\phi g'=0$ and $\theta f' = g' p$,
so if and only if $g'=\theta h$ and $f'=h p$
for some $h\in\Hom(V,M)$.
Thus $\pi^{-1}(x,\theta,\phi)$ is
isomorphic to $\Hom(V,M)$.
It follows that 
$\dim \mathcal{D} = \dim\mathcal{E}+\dim\Hom(V,M)$.

There is also a projection 
$\mathcal{D}\to\Hom(U,M)$. 
To compute the fibre 
over $f\in\Hom(U,M)$, we construct the pushout
\[
\begin{CD}
0 @>>> U @>p>> V @>q>> L @>>> 0 \\
@. @VfVV @VaVV @\vert @. \\
0 @>>> M @>r>> N @>s>> L @>>> 0.
\end{CD}
\]
Any invertible linear map $b:N\to K^\alpha$
induces an $A$-module structure on $K^\alpha$
such that $b$ is a module homomorphism.
If $b$ is given by invertible linear maps $e_i N\to K^{\alpha_i}$
($i\in I$), this structure corresponds
to an element $x\in\Rep(A,\alpha)$ with
$b\in\Hom(N,K_x)$ and 
$((x,br,sb^{-1}),f,ba)\in\mathcal{D}$.
Conversely, if $((x,\theta,\phi),f,g)$ is in 
$\mathcal{D}$, then by the universal property of the pushout 
we have constructed, there 
is a map $b\in\Hom(N,K_x)$ giving a commutative 
diagram
\[
\begin{CD}
0 @>>> M @>r>> N @>s>> L @>>> 0 \\
@. @\vert @VbVV @\vert @. \\
0 @>>> M @>\theta>> K_x @>\phi>> L @>>> 0.
\end{CD}
\]
and clearly $b$ is an isomorphism. Thus the 
fibre over $f$ is isomorphic to 
the product over $i\in I$ of the set of invertible
linear maps $e_i N\to K^{\alpha_i}$. 
Thus each fibre is irreducible of dimension $\alpha\cdot\alpha$,
so 
\[
\dim\mathcal{D} = \dim\Hom(U,M) + \alpha\cdot\alpha.
\]
The assertion follows.
\end{proof}

\begin{lemma}\label{l:onetopdim}
Let $\mathcal{C}$ be a $\GL(\gamma)$-stable constructible
subset of $\Rep(A,\gamma)$, let $T$ be a simple
$A$-module of dimension vector $\beta$, let $m$
be a positive integer, and let $\alpha=\gamma+m\beta$.
Let $\mathcal{F}$ be the set of $x\in\Rep(A,\alpha)$
such that $\dim\Hom(K_x,T)=m$ and there is an
exact sequence
\[
0\to K_c\xrightarrow{\theta} K_x\xrightarrow{\phi} T^m\to 0
\]
with $c\in \mathcal{C}$. Then $\mathcal{F}$ is constructible and
\[
\dim\mathcal{F} \le \dim\mathcal{C}
+ \alpha\cdot\alpha
- \gamma\cdot\gamma
- m^2 + md
\]
where
\[
d = \max_{c\in \mathcal{C}} \left( \dim\Ext^1(T,K_c) - \dim\Hom(T,K_c)\right).
\]
\end{lemma}

\begin{proof}
Let $\mathcal{G}$ be the set of quadruples 
\[
(x,c,\theta,\phi) \in \Rep(A,\alpha)\times\Rep(A,\gamma)\times\Hom_K(K^\gamma,K^\alpha)
\times\Hom_K(K^\alpha,T^m)
\]
forming an exact sequence as above, with $c\in\mathcal{C}$, 
and with $\dim\Hom(K_x,T)=m$.
Clearly $\mathcal{G}$ is constructible. The first 
projection $\mathcal{G}\to\Rep(A,\alpha)$ 
has image $\mathcal{F}$, so that $\mathcal{F}$ is constructible.
If $(x,c,\theta,\phi)\in\mathcal{G}$, since the 
components of $\phi$ must form a basis of $\Hom(K_x,T)$,
any other surjective map $K_x\to T^m$ must be obtained
by composing $\phi$ with an automorphism of $T^m$. 
It follows that $\GL(\gamma)\times\Aut_A(T^m)$ acts freely
and transitively on the fibres of the projection, so
\[
\dim \mathcal{F} = \dim\mathcal{G} - \gamma\cdot\gamma - m^2.
\]
The second projection $\mathcal{G}\to\Rep(A,\gamma)$
has image contained in $\mathcal{C}$, and the 
fibre over any point $c$ is contained in the variety
$\mathcal{E}$ of Lemma \ref{l:extE} for $L=T^m$ and
$M=K_c$. It follows that the fibre has dimension at most
$\alpha\cdot\alpha+md$. Thus 
$\dim\mathcal{G}\le\dim\mathcal{C}+\alpha\cdot\alpha+md$,
giving the result.
\end{proof}

\section{Dimensions of strata}
\label{s:dimstrata}

Let $\Pi = \Pi^0$ be the (undeformed) preprojective algebra for $Q$. 
The trivial paths in $K\overline{Q}$
give a complete set of orthogonal idempotents $e_i$
in $\Pi$, and for $\alpha\in\N^I$ we can identify $\Rep(\Pi,\alpha)$
with $\mu_\alpha^{-1}(0)$.

Let $T_1,\dots,T_r$ be a collection of non-isomorphic simple $\Pi$-modules,
of dimensions $\beta^{(1)},\dots\beta^{(r)}$.

If $h\ge 0$, $j_1,\dots,j_h$ are integers in the range $1,\dots,r$,
and $m_1,\dots,m_h$ are positive integers, we
say that a $\Pi$-module $M$ has \emph{top-type}
$(j_1,m_1; \dots; j_h,m_h)$ if it has a
filtration by submodules
\[
0 = M_0 \subset M_1 \subset \dots \subset M_h = M
\]
such that $M_i/M_{i-1} \cong T_{j_i}^{m_i}$
and $\dim\Hom(M_i,T_{j_i}) = m_i$ for each $i$.
This means that having chosen the sequence $j_1,\dots,j_h$,
the submodules $M_i$ are determined by decreasing induction
on $i$ by the fact that $M_{i-1}$ is the intersection
of the kernels of all maps from $M_i$ to $T_{j_i}$ (the `top' of $M_i$ with
respect to the simple $T_{j_i}$).

If $(j_1,m_1; \dots; j_h,m_h)$
is a top-type and $1\le i\le h$, define 
$z_i$ to be zero if $p(\beta^{(j_i)})=0$
or there is no $k<i$ with $j_k=j_i$, and otherwise 
to be equal to $m_k$ for the largest $k<i$ with $j_k=j_i$.

\begin{lemma}\label{l:zbounds}
If $M$ has top-type
$(j_1,m_1; \dots; j_h,m_h)$ then
$\dim\Hom(M_{i-1},T_{j_i})\le z_i$.
\end{lemma}

\begin{proof}
For simplicity of notation, let $T=T_{j_i}$.
First suppose that $\Ext^1(T,T)=0$.
Applying $\Hom(-,T)$ to the exact sequence
$0\to M_{i-1}\to M_i\to T^{m_i}\to 0$
gives an exact sequence
\[
0\to\Hom(T^{m_i},T)\to\Hom(M_i,T)\to\Hom(M_{i-1},T)
\to\Ext^1(T^{m_i},T) = 0.
\]
Thus $\Hom(M_{i-1},T)=0$ by dimensions.

Now suppose that $\Ext^1(T,T)\neq 0$.
If there is no $k<i$ with $j_k=j_i$ then $M_{i-1}$
has no composition factors isomorphic to $T$,
so the Hom space is zero. Thus suppose otherwise,
and choose $k<i$ maximal with $j_k=j_i$.
By the definition of top-type we have
$\dim\Hom(M_k,T) = m_k$. Now
applying $\Hom(-,T)$ to the exact sequence
\[
0\to M_k\to M_{i-1}\to M_{i-1}/M_k\to 0
\]
and using the fact that $M_{i-1}/M_k$ has no composition
factor isomorphic to $T$, we deduce that
$\Hom(M_{i-1},T)$ embeds in $\Hom(M_k,T)$.
\end{proof}

\begin{lemma}
The subset $\mathcal{D}$ 
of $\mu_\alpha^{-1}(0)$ consisting of those 
elements with top-type $(j_1,m_1; \dots; j_h,m_h)$
is constructible and has dimension at most
\[
\alpha\cdot\alpha-1+p(\alpha)
+ \sum_{i=1}^h m_i z_i
- \sum_{i=1}^h m_i^2 p(\beta^{(j_i)}).
\]
\end{lemma}

\begin{proof}
Let $\gamma=\alpha-m_h\beta^{(j_h)}$.
Let $\mathcal{C}$ be the subset of 
$\mu_\gamma^{-1}(0)=\Rep(\Pi,\gamma)$
consisting of those elements $c$
such that $K_c$ has top-type
$(j_1,m_1; \dots; j_{h-1},m_{h-1})$
and with $\dim\Hom(K_c,T_{j_h})\le z_h$.
By induction on $h$ it is constructible and 
\[
\dim\mathcal{C} \le 
\gamma\cdot\gamma-1+p(\gamma)
+ \sum_{i=1}^{h-1} m_i z_i
- \sum_{i=1}^{h-1} m_i^2 p(\beta^{(j_i)}).
\]

By Lemma~\ref{l:zbounds}, for any 
element $x\in\mathcal{D}$ there is an
exact sequence
\[
0\to K_c\to K_x\to T_{j_h}^{m_h}\to 0
\]
with $c\in\mathcal{C}$. Thus by
Lemma~\ref{l:onetopdim} we have
\[
\dim\mathcal{D} \le \dim\mathcal{C}
+ \alpha\cdot\alpha
- \gamma\cdot\gamma
- m_h^2 + m_h d
\]
where
\[
d = \max_{c\in \mathcal{C}} \left( \dim\Ext^1(T_{j_h},K_c) 
- \dim\Hom(T_{j_h},K_c)\right).
\]
Now by \cite{CBex} we have
\[
\dim\Ext^1(T_{j_h},K_c) - \dim\Hom(T_{j_h},K_c)
= \dim\Hom(K_c,T_{j_h}) - (\gamma,\beta^{(j_h)}),
\]
so $d\le z_h - (\gamma,\beta^{(j_h)})$.
The result follows.
\end{proof}

Let $\pi:\mu_\alpha^{-1}(0)\to N_Q(0,\alpha)$
be the quotient map.
Fix a representation type
$(k_1,\beta^{(1)};\dots ;k_r,\beta^{(r)})$
for representations
of $\Pi$ of dimension $\alpha$.

\begin{theorem}
If $x\in N_Q(0,\alpha)$ has representation type 
$(k_1,\beta^{(1)};\dots ;k_r,\beta^{(r)})$,
then
\[
\dim \pi^{-1}(x)\le 
\alpha\cdot\alpha - 1 +p(\alpha) - \sum_t p(\beta^{(t)}).
\]
\end{theorem}

\begin{proof}
Any element of $\pi^{-1}(x)$ has top-type
$(j_1,m_1; \dots; j_h,m_h)$
for suitable $j_i$ and $m_i$.
It thus suffices to prove that
\[
\sum_{i=1}^h \left( m_i z_i - m_i^2 p(\beta^{(j_i)}) \right) 
\le -\sum_{t=1}^r p(\beta^{(t)}).
\]
for each possible top-type.
We group the first sum according to the value of $j_i$.
Fix $t$, and let the values of $i$ with $j_i=t$ be
$i_1,\dots,i_q$. It clearly suffices to prove that
\[
\sum_{p=1}^q 
\left( m_{i_p} z_{i_p} - m_{i_p}^2 p(\beta^{(t)}) \right) 
\le -p(\beta^{(t)}).
\]

If $p(\beta^{(t)})=0$, then $z_{i_p} = 0$,
and the inequality is trivial. Thus suppose
that $p(\beta^{(t)}) = \ell > 0$.
Defining $n_p = m_{i_p}$ for $1\le p\le q$,
we have $z_{i_p} = n_{p-1}$ if $p>1$, and $z_{i_1}=0$.
Thus the inequality becomes
\[
\sum_{p=2}^q n_p n_{p-1} - \sum_{p=1}^q n_p^2 \ell \le -\ell
\]
or equivalently
\[
(\ell-1)(\sum_{p=1}^q n_p^2) +
\left( \sum_{p=1}^q n_p^2 - \sum_{p=2}^q n_p n_{p-1} \right)
\ge \ell.
\]
Now the first of the two terms on the left 
is clearly at least $\ell-1$,
and the second is at least 1, for it is the (positive
definite) quadratic form for the Dynkin graph of type $\mathbb{A}_q$ 
applied to the dimension vector $(n_1,\dots,n_q)$.
\end{proof}

In view of \cite[Theorem 1.3]{CBmm}, this has an immediate corollary.

\begin{corollary}
\label{c:pistratumdim}
The inverse image in $\mu_\alpha^{-1}(0)$
of the stratum of representation type 
$(k_1,\beta^{(1)};\dots ;k_r,\beta^{(r)})$
has dimension at most
$\alpha\cdot\alpha - 1 +p(\alpha) + 
\sum_t p(\beta^{(t)})$.
\end{corollary}

If $Q$ has no oriented cycles and the $\beta^{(t)}$
are the coordinate vectors, then the stratum
consists only of the origin, and its inverse image in
$\mu_\alpha^{-1}(0)$ is the nilpotent cone.
According to Lusztig \cite{Lusztig}
it is Lagrangian, so its dimension is equal to the bound 
in the corollary.

\section{Normality of quotients}
\label{s:normquot}

Let $X$ be an irreducible affine scheme of finite type over $K$, let $G$ be a 
reductive group acting on $X$, and let $\pi:X\to X \dbslash G$ be the natural map.
It is well known that if $X$ is normal, then so is $X \dbslash G$. 
Here we find a different condition which still implies normality of the quotient.

Recall \cite[\S 5.7]{EGAIVtwo} that a scheme $X$ 
is said to have property $(S_i)$ if
\[
\depth(\mathcal{O}_{X,Z}) \ge \min(i,\codim_X Z)
\]
for every irreducible closed subset $Z$ of $X$.
By \cite[Theorem 3.8]{Localcohom}, for affine $X$ it is
equivalent that the local cohomology groups $H_Z^j(X)$
vanish for all such $Z$ and all $j < \min(i,\codim_X Z)$.

\begin{theorem}
Suppose that $X$ and an open subscheme $U$ of $X \dbslash G$
have property $(S_i)$. If $\pi^{-1}(X \dbslash G \setminus U)$ has
codimension at least $i$ in $X$, then $X \dbslash G$ has $(S_i)$. 
\end{theorem}

\begin{proof}
Suppose that $Z$ is an irreducible closed
subset of $X \dbslash G$. Let $V$ be the complement of $U$.
By \cite[Proposition 1.9]{Localcohom} there is a long exact
sequence
\[
0 \to H_{Z\cap V}^0 (X \dbslash G) 
\to H_Z^0 (X \dbslash G) 
\to H_{Z\cap U}^0 (X \dbslash G)
\to H_{Z\cap V}^1 (X \dbslash G) 
\to \dots,
\]
so it suffices to prove that if
$j < \min(i,\codim_{X \dbslash G} Z)$
then 
$H_{Z\cap V}^j (X \dbslash G)$ and 
$H_{Z\cap U}^j (X \dbslash G)$
both vanish.

Now $H_{Z\cap V}^j (X \dbslash G)$
embeds in $H_{\pi^{-1}(Z\cap V)}^j (X)$
by \cite[Corollary 6.8]{HochRob},
and this space vanishes since $X$ has property $(S_i)$
and $\pi^{-1}(Z\cap V)$ has codimension at least $i$ in $X$.

Also $H_{Z\cap U}^j (X \dbslash G) \cong H_{Z\cap U}^j (U)$
by excision \cite[Proposition 1.3]{Localcohom},
and this vanishes since $U$ has $(S_i)$.
\end{proof}

Recall that a scheme has property $(R_i)$ if its singular locus 
has codimension at least $i$. Serre's criterion says that 
$X$ is normal if and only if it has $(R_1)$ and $(S_2)$.

\begin{corollary}
\label{c:quotnormal}
Let $U$ be an open subscheme of $X \dbslash G$. 
If $U$ is normal, its complement $Z$ has codimension at
least two in $X \dbslash G$, $\pi^{-1}(Z)$ has
codimension at least two in $X$, and $X$ has $(S_2)$,
then $X \dbslash G$ is normal.
\end{corollary}

\begin{proof}
$X \dbslash G$ is regular in codimension 1 since $U$ 
is normal and and its complement has codimension at least 2. 
Thus we just need $X \dbslash G$ to have property $(S_2)$,
and this is a special case of the theorem.
\end{proof}

\section{Completion of the proof}
\label{s:completion}

Recall that a dimension vector $\alpha$ is \emph{sincere} if $\alpha_i>0$ for all $i$.
Given a quiver $Q$ and a sincere dimension vector $\alpha$, we say
that we are in the \emph{nearly Kleinian} case if either

(1) $Q$ has only one vertex, say $i$, and $\alpha_i=1$, or

(2) $Q$ is obtained from an extended Dynkin quiver by adjoining any number
of loops at extending vertices, and $\alpha = \delta$, the minimal positive
imaginary root for the extended Dynkin quiver. (Recall that $i$ is an
\emph{extending vertex} if $\delta_i = 1$.)

\begin{lemma}
$N_Q(0,\alpha)$ is normal in the nearly Kleinian case.
\end{lemma}

\begin{proof}
Observe that the preprojective relations involve commutators, and if $a$ is a loop
and the dimension vector at the corresponding vertex is 1, then $aa^*-a^*a=0$
automatically. Thus $a$ and $a^*$ can be removed from the relations.
The result is now clear---in case (1) $N_Q(0,\alpha)$ is an affine space,
and in case (2) it is the product of a Kleinian singularity and an affine space.
\end{proof}

\begin{lemma}
Let $\alpha$ be a sincere dimension vector for a quiver $Q$ with
$\alpha\in\Sigma_0$, and suppose that we are not in the nearly Kleinian case.
If $(k_1,\beta^{(1)}; \dots ; k_r,\beta^{(r)})$ is a representation type
for dimension vector $\alpha$, then at least one of the following holds.

(i) $\sum_t k_t^2 < \sum_i \alpha_i^2$.

(ii) $p(\alpha) - \sum_t p(\beta^{(t)}) \ge 2$.
\end{lemma}

\begin{proof}
We suppose that (i) and (ii) fail, and derive a contradiction.
We have
$\sum_i \alpha_i^2 = \alpha\cdot\alpha = \sum_{s,t} k_s k_t \beta^{(s)}\cdot\beta^{(t)} \ge \sum_t k_t^2$.
Since (i) fails this must be an equality, which since the $\beta^{(t)}$ have non-negative integer components,
implies that the $\beta^{(t)}$ must be distinct coordinate vectors. 
Since $\alpha$ is sincere, they
are exactly the coordinate vectors in some order, and the multiplicities $k_t$
are  the corresponding components of $\alpha$.

If $\alpha=\beta^{(t)}$ for some $t$ then we are in the nearly Kleinian case, a contradiction.
Thus the decomposition of $\alpha$ as a sum over
$t$ of $k_t$ copies of each $\beta^{(t)}$ has at least two terms, 
so since $\alpha\in\Sigma_0$ we know that 
$p(\alpha) > \sum_t k_t p(\beta^{(t)})$. Since (ii) fails, it follows that
$p(\alpha) = 1 + \sum_t k_t p(\beta^{(t)}) = 1 + \sum_t p(\beta^{(t)})$, and
$k_t=1$ whenever $p(\beta^{(t)}) > 0$, that is,
whenever there is a loop at the corresponding vertex.

Now we have
\begin{align*}
q(\alpha) - \sum_i k_t q(\beta^{(t)})
&= \sum_t q(k_t \beta^{(t)}) + \frac12 \sum_{\substack{s,t \\ s\neq t}} (k_s \beta^{(s)},k_t \beta^{(t)}) - \sum_t k_t q(\beta^{(t)})
\\
&= \sum_t (k_t^2 - k_t) q(\beta^{(t)}) + \frac12 \sum_{\substack{s,t \\ s\neq t}} (k_s \beta^{(s)},k_t \beta^{(t)})
\\
&= \frac12 \sum_t k_t \biggl((k_t-1)(\beta^{(t)},\beta^{(t)}) + \sum_{\substack{s \\ s\neq t}} (k_s \beta^{(s)},k_t \beta^{(t)})\biggr)
\\
&= \frac12 \sum_t k_t (\beta^{(t)},\alpha-\beta^{(t)}).
\end{align*}
As observed above, we have $\beta^{(t)}\neq \alpha$, so
$(\beta^{(t)},\alpha-\beta^{(t)}) = -2-d_t$ for some $d_t\ge 0$ by \cite[Corollary 5.7]{CBmm}.
Then 
\[
p(\alpha) - \sum_t k_t p(\beta^{(t)})
= 1 - \sum_t k_t + \frac12 \sum_t k_t(2+d_t) = 1 + \frac12 \sum_t k_t d_t.
\]
Now since (ii) fails we have $d_t = 0$ for all $t$, that is, that 
$(\beta^{(t)},\alpha-\beta^{(t)}) = - 2$.
Equivalently
\[
(k_t-1)(\beta^{(t)},\beta^{(t)}) + \sum_{\substack{s \\ s\neq t}} k_s (\beta^{(s)},\beta^{(t)}) = -2.
\]

Let $Q''$ be the quiver obtained from $Q$ by deleting all loops.
Since the $\beta^{(t)}$ are distinct coordinate vectors we have
\begin{align*}
(\beta^{(t)},\alpha)_{Q''} 
&= k_t (\beta^{(t)},\beta^{(t)})_{Q''} + \sum_{\substack{s \\ s\neq t}} k_s (\beta^{(t)},\beta^{(s)})_{Q''}
\\
&= 2 k_t + \sum_{\substack{s \\ s\neq t}} k_s (\beta^{(t)},\beta^{(s)})_Q
= 2 k_t - 2 - (k_t-1)(\beta^{(t)},\beta^{(t)})_Q.
\end{align*}
Now this is zero, for either $(\beta^{(t)},\beta^{(t)})_Q=2$, or there are loops in $Q$ at
the vertex, in which case $k_t=1$.

Since $\alpha$ is a vector for $Q''$ with positive integer components, 
connected support, and the
bilinear form with each coordinate vector gives 0, it follows that $Q''$ is
an extended Dynkin quiver and $\alpha$ is a multiple of the minimal
positive imaginary root $\delta$ (see \cite[\S 1.2]{Kac}). Now if $Q$ contains any loops
then some $k_t=1$, so that $\alpha=\delta$. Otherwise $Q=Q''$ is
extended Dynkin itself,
and the fact that $\alpha\in\Sigma_0$ implies that $\alpha=\delta$.
Thus we are in the nearly Kleinian case, a contradiction.
\end{proof}

We can now prove our main result.

\begin{theorem}
$N_Q(\lambda,\alpha)$ is a normal variety for any
$\lambda$ and $\alpha$.
\end{theorem}

\begin{proof}
In view of the etale local structure result,
it suffices to prove normality for $\lambda=0$.
We prove it for all quivers and dimension vectors $\alpha$ 
by induction on $\sum_i \alpha_i^2$.
Since a product of normal varieties is normal, and a
symmetric product of copies of a normal variety is normal,
by the main result of \cite{CBde} we may assume that $\alpha\in\Sigma_0$.
Clearly we may suppose that $\alpha$ is sincere.
Also, we may assume that we are not in the nearly Kleinian case,
since that has been dealt with separately.

By the lemma, for any
representation type $(k_1,\beta^{(1)}; \dots ; k_r,\beta^{(r)})$, 
conditions (i) or (ii) of the lemma hold.

By the etale local structure, a point in the stratum
of this representation type has a neighbourhood which is isomorphic 
in the etale topology to a neighbourhood of 0 in $N_{Q'}(0,\kappa)$.
By the remark before Theorem \ref{th:ettopcoordfree},
the points of strata where (i) holds are normal by induction.

On the other hand, if a stratum satisfies (ii), then
Corollary \ref{c:pistratumdim} and \cite[Theorem 1.2]{CBmm}
imply that the inverse image of the stratum in 
$\mu_\alpha^{-1}(0)$ has 
codimension at least 2.

We now apply Corollary \ref{c:quotnormal}.
Let  $U$ be the normal locus of $N_Q(0,\alpha)$, an open subset.
By the observations above, its complement has inverse
image in $\mu_\alpha^{-1}(0)$ of codimension at least 2.
Thus in fact $N_Q(0,\alpha)$ is normal.
\end{proof}

Using Geometric Invariant Theory quotients it is 
possible to construct resolutions of singularities. 
If $\theta\in \Z^I$ satisfies $\theta\cdot\alpha=0$, then there are notions
of $\theta$-stable and $\theta$-semistable elements of $\mu_\alpha^{-1}(\lambda)$,
there is a GIT quotient $\mu_\alpha^{-1}(\lambda) \dbslash (\GL(\alpha),\theta)$,
and a semisimplification map
\[
\mu_\alpha^{-1}(\lambda) \dbslash (\GL(\alpha),\theta) \to N_Q(\lambda,\alpha)
\]
which is a projective morphism, see \cite{King}. Moreover we have:

(1) If $\alpha\in\Sigma_\lambda$, then $\mu_\alpha^{-1}(\lambda)$ 
is irreducible by \cite{CBmm}, and the general element is a simple 
representation of $\pr$, hence $\theta$-stable. Thus the semisimplification 
map is a birational morphism of irreducible varieties.

(2) If $\alpha$ is indivisible, that is, its components have no common divisor, 
then one can choose a weight $\theta$ such that every $\theta$-semistable element 
is $\theta$-stable. Since any $\theta$-stable element has trivial endomorphism
algebra as a $\pr$-module, the fibre $\mu_\alpha^{-1}(\lambda)$ is smooth
(see for example \cite[Lemma 10.3]{CB}), and then $N_Q(\lambda,\alpha)$ is
smooth by Luna's Slice Theorem \cite[Corollary 1, p98]{Luna}.

Thus, if (1) and (2) hold, then the semisimplification map is
a resolution of singularities. 
As usual, the Grauert-Riemenschneider Theorem 
gives the following consequence (cf.\ \cite[p50]{TE}).

\begin{corollary}
If $\alpha\in\Sigma_\lambda$ is indivisible, then $N_Q(\lambda,\alpha)$
has rational singularities.
\end{corollary}

\begin{proof}
It suffices to observe that $N_Q(\lambda,\alpha)$ is normal and that
if $\theta$ is chosen so that $\theta$-semistables are $\theta$-stable, then since
$\mu_\alpha^{-1}(\lambda) \dbslash (\GL(\alpha),\theta)$ 
is a symplectic quotient, it has a symplectic form, and hence has trivial canonical class. 
\end{proof}

Incidentally, if $\alpha\in\Sigma_0$, then the moment map 
$\mu_\alpha$ is flat by \cite{CBmm}, and $\alpha\in\Sigma_\lambda$ 
for any $\lambda$ with $\lambda\cdot\alpha=0$. 
Letting $\mathfrak{h} = \{\lambda\in K^I \mid \lambda\cdot\alpha=0\}$,
there is a flat family $\mu_\alpha^{-1}(\mathfrak{h})\dbslash \GL(\alpha)\to \mathfrak{h}$
whose fibre over 0 is $N_Q(0,\alpha)$.
In case $\alpha$ is also indivisible, the general fibre of this family is smooth, and 
the whole family has a simultaneous resolution of singularities 
$\mu_\alpha^{-1}(\mathfrak{h})\dbslash (\GL(\alpha),\theta)$.
If $Q$ is extended Dynkin and $\alpha=\delta$, then the family is
the semi-universal deformation, or rather its lift through a Weyl group action.

One may wonder whether $N_Q(0,\alpha)$ is actually singular.
We have a partial answer based on the following lemma.

\begin{lemma}
Suppose that $\pi:\tilde{X}\to X$ is a projective, birational morphism between smooth
irreducible varieties. If $\tilde{X}$ has trivial canonical class, and
$\pi^{-1}(U)\to U$ is an isomorphism for some open $U\subseteq X$ whose
complement $Z$ has codimension $\ge 2$, then $\pi$ is an isomorphism.
\end{lemma}

\begin{proof}
Let $j:U\to X$ be the inclusion, and let $i:U\to \tilde{X}$ be the composition of the isomorphism
$U\to \pi^{-1}(U)$ and the inclusion $\pi^{-1}(U)\to\tilde{X}$.

Since $X$ is smooth, the codimension condition implies that any locally 
free sheaf $\mathcal{F}$ on $X$ is $Z$-closed \cite[Theorem 5.10.5]{EGAIVtwo}, 
so the map $j^* : \Gamma(X,\mathcal{F}) \to  \Gamma(U,\mathcal{F}|_U)$
is an isomorphism by \cite[Proposition 5.9.8]{EGAIVtwo}. Similarly, if $\mathcal{G}$
is a locally free sheaf on $\tilde{X}$ then
$i^* : \Gamma(\tilde{X},\mathcal{G}) \to  \Gamma(U,\mathcal{G}|_U)$
is injective by $Z$-purity \cite[Proposition 5.10.2]{EGAIVtwo}.

Let $X$ and $\tilde{X}$ have dimension $n$. Since $\tilde{X}$ has trivial canonical 
class, there is a nowhere vanishing $n$-form $\sigma$ on $\tilde{X}$.
Then $i^*\sigma$ is such a form on $U$, and since 
$j^* : \Gamma(X,\omega_X) \to  \Gamma(U,\omega_U)$ is an isomorphism,
there is an $n$-form $\tau$ on $X$ with $j^*\tau = i^*\sigma$. 
Now $i^* \pi^* \tau = j^* \tau = i^* \sigma$, so the fact that
$i^* : \Gamma(\tilde{X},\omega_{\tilde{X}}) \to  \Gamma(U,\omega_U)$
is injective implies that $\pi^* \tau = \sigma$.

In particular $\pi^*\tau$ is a nowhere vanishing $n$-form on $\tilde{X}$. 
Now if $x\in\tilde{X}$ then 
$(\pi^*\tau)_x$ is the composition
\[
\wedge^n T_x \tilde{X} \xrightarrow{\wedge^n d\pi_x} \wedge^n T_{\pi(x)} X \xrightarrow{\tau_{\pi(x)}} K
\]
where $d\pi_x : T_x \tilde{X} \to T_{\pi(x)} X$ is the induced map on tangent spaces.
Since this is nonzero, $d\pi_x$ must be invertible. Thus $\pi$ is etale, and hence it
has finite fibres. Then by Zariski's Main Theorem, $\pi$ is bijective. Since it is
also birational, it is an isomorphism.
\end{proof}

\begin{proposition}
Let $Q$ have no loops at vertices, and let $\alpha\in\Sigma_0$ be indivisible.
If $\alpha$ is an imaginary root (equivalently, not a coordinate vector),
then $N_Q(0,\alpha)$ is singular.
\end{proposition}

\begin{proof}
Since $N_Q(0,\alpha)$ does not depend on the orientation of $Q$, 
see for example \cite[Lemma 2.2]{CBH}, and since $Q$ has no loops, we can reorient
it so that it has no oriented cycles.

Choosing $\theta$ so that any $\theta$-semistable is $\theta$-stable, 
the semisimplification map is a resolution of singularities, and
the quotient $\mu_\alpha^{-1}(\lambda) \dbslash (\GL(\alpha),\theta)$
classifies $\GL(\alpha)$-orbits of $\theta$-stable elements of $\mu_\alpha^{-1}(0)$. 
By the lemma, if $N_Q(0,\alpha)$ is smooth, 
the semisimplification map is a bijection. Thus, considering the fibre
over 0, the set of nilpotent $\theta$-stable elements of $\mu_\alpha^{-1}(0)$
must form a single orbit, so must have dimension $\dim \GL(\alpha) - 1$ 
(taking into account the fact
that $\GL(\alpha)$ has a subgroup isomorphic to $K^*$ which acts trivially).

Since $Q$ has no oriented cycles, by \cite{Lusztig} the nilpotent elements of 
$\mu_\alpha^{-1}(0)$ form an equidimensional subvariety of 
$\Rep(\overline{Q},\alpha)$ of dimension $\dim \Rep(Q,\alpha)$. Since 
the nilpotent elements which are $\theta$-stable
form an open subset of this, we must have 
\[
\dim \GL(\alpha)-1 = \dim\Rep(Q,\alpha).
\]
In other words $q(\alpha)=1$, so that $\alpha$ is a real root. The equivalence with
coordinate vectors follows from the results in \cite[Section 7]{CBmm}.
\end{proof}

\section{Appendix}

We show that Kraft and Procesi's construction of
conjugacy class closures \cite{KP} works in more than just the nilpotent case.
Let $C$ be a conjugacy class in $\Mat(n,K)$, and choose
$t\ge 1$ and elements $\xi_1,\dots,\xi_t\in K$ such that
$\prod_{j=1}^t (M-\xi_j 1)=0$ for $M\in C$.
For $i\le t$ let $n_i$ be the rank of $\prod_{j=i+1}^t (M-\xi_j 1)$,
so $n_t = n$ and $n_0 = 0$.
Define $Z$ as in the introduction, but with the relations
$B_1 A_1 = (\xi_1 - \xi_2) 1$ and $B_i A_i = A_{i-1} B_{i-1} + (\xi_i - \xi_{i+1}) 1$
for $2\le i\le t-1$. Let $H$ be as in the introduction.

\begin{lemma}
\label{l:kpgen}
The map $\Theta:Z\to\Mat(n,K)$ sending $(A_i,B_i)$ to $A_{t-1}B_{t-1}+\xi_t 1$
induces an isomorphism
$Z \dbslash H \cong \overline{C}$.
\end{lemma}

\begin{proof}
Clearly $\Theta$ is constant on orbits of $H$ and equivariant for the
conjugation action of $\GL(n,K)$. The First Fundamental Theorem of 
Invariant Theory and an induction on $t$ shows that $\Theta$ induces a 
closed embedding $Z\dbslash H\to\Mat(n,K)$, see \cite[\S 2]{KP}.
The image of $\Theta$ contains $C$, since if $M\in C$ one can
identify $K^{n_i}$ with the image of $\prod_{j=i+1}^t (M-\xi_j 1)$,
and take $B_i = M-\xi_{i+1}1$ and $A_i$ to be the inclusion.
In \cite{KP} the proof is completed using the Gerstenhaber-Hesselink
Theorem, and perhaps our lemma could be proved in a similar way; here we use 
the Marsden-Weinstein reductions. 
Clearly it suffices to prove that $Z \dbslash H$ is irreducible of dimension
\[
\dim C = n^2 - \sum_{i=1}^t d_i^2
\]
where $d_i = n_i - n_{i-1}$.
By deframing, $Z\dbslash H \cong N_Q (\lambda,\alpha)$ for 
the quiver $Q$ illustrated in the introduction, 
$\alpha=(n_1,\dots,n_{t-1},1)$
and $\lambda=(\xi_2-\xi_1,\dots,\xi_t-\xi_{t-1},\nu)$,
where $\nu$ is chosen so that $\lambda\cdot\alpha=0$.
Clearly
\[
p(\alpha) = \sum_{i=1}^{t-1} (n_i n_{i+1} - n_i^2) = \frac12 (n^2 - \sum_{i=1}^t d_i^2),
\]
so by \cite[Corollary 1.4]{CBmm} it suffices to show that $\alpha\in\Sigma_\lambda$.

Any decomposition of $\alpha$ as a sum of positive roots must include one 
term whose component at vertex $t$ is nonzero, say $\beta = (m_1,\dots,m_{t-1},1)$
with $0\le m_1\le \dots\le m_{t-1}\le n$.
All other terms must be of the form $\gamma^{(jk)}$ with $1\le j\le k\le t-1$, where
\[
\gamma^{(jk)}_i =  
\begin{cases} 1 & \text{(if $j\le i\le k$)} \\ 0 & \text{(else).}
\end{cases}
\]
The decompositions involved in the definition of $\Sigma_\lambda$ are 
those in which the summands satisfy $\lambda\cdot\gamma^{(jk)}=0$,
that is, $\xi_j = \xi_{k+1}$. 

Let $e_i = m_i-m_{i-1}$ (with $m_t = n$ and $m_0=0$).
Clearly the sequence $e_1,\dots,e_t$ is obtained from $d_1,\dots,d_t$
by decreasing $d_j$ by one, and increasing $d_{k+1}$ by one, for each
term of the form $\gamma^{(jk)}$ in the decomposition of $\alpha$.

Since $p(\gamma^{(jk)}) = 0$, to show that $\alpha\in \Sigma_\lambda$
we need $p(\alpha) > p(\beta)$ for a nontrivial decomposition,
or equivalently $\sum_{i=1}^t d_i^2 < \sum_{i=1}^t e_i^2$.
By including additional terms of the form $\gamma^{(jk)}$ in the decomposition, 
we may reorder the $e_i$, and hence assume that 
$e_j \le e_{\ell}$ whenever $\xi_j = \xi_{\ell}$. This will reduce the
size of $\beta$, but not change the value of $\sum_{i=1}^t e_i^2$.
We now partition the sums according to the value of $\xi_i$. 
If $\mu\in K$ and $\{i \mid \xi_i = \mu\} = \{i_1,i_2,\dots,i_p\}$ with 
$i_1<i_2<\dots<i_p$, then $d_{i_{p-q}}$ is the number of Jordan blocks
of eigenvalue $\mu$ of size $> q$, so $d_{i_1} \le d_{i_2} \le \dots \le d_{i_p}$.
It follows that $\sum_{j=1}^p d_{i_j}^2 \le \sum_{j=1}^p e_{i_j}^2$, 
strict if not all $e_{i_j} = d_{i_j}$. This is just a numerical condition,
but to avoid writing out a proof, observe that it is implicit in the
Gerstenhaber-Hesselink Theorem (see \cite[\S 1.3]{KP}). 
\end{proof}

\begin{remark}
\label{r:notnew}
In the definition $N_Q(\lambda,\alpha) = \mu_\alpha^{-1}(\lambda) \dbslash \GL(\alpha)$
we have only considered Marsden-Weinstein reductions associated to fixed points 
of the dual of the Lie algebra. One could more generally consider 
$\mu_\alpha^{-1}(\overline{C}) \dbslash \GL(\alpha)$ 
for an arbitrary coadjoint orbit $C$. Here, by taking the closure of the coadjoint orbit, 
the inverse image is an affine variety, 
so the quotient is well-defined.
However, the coadjoint orbit is determined by fixing a conjugacy class in
$\Mat(\alpha_i,K)$ for each vertex $i$, and by using Lemma \ref{l:kpgen} one sees that
the more general Marsden-Weinstein reduction is isomorphic 
$N_{Q'}(\lambda',\alpha')$ for a quiver $Q'$ which is obtained from $Q$ by
attaching a linear quiver to each vertex of $Q$. 

The same comment applies to the Marsden-Weinstein
reduction $N(A,M,\omega,\zeta)$ considered in Section \ref{s:coordfree}.

Actually this construction is already
implicit in \cite[\S 1]{W}, where $Q$ consists of one vertex and a single loop, and
$Q'$ is obtained by attaching a new vertex to the original one by a single edge.
\end{remark}

\end{document}